\documentstyle[12pt]{article}
\newlength{\abstractwidth}
\renewcommand{\thefootnote}{\fnsymbol{footnote}}
\renewcommand{\thanks}[1]{\footnote{#1}} 
\newcommand{\starttext}{ \setcounter{footnote}{0}
\renewcommand{\thefootnote}{\arabic{footnote}}}

\newcommand{\be}{\begin{equation}}
\newcommand{\bea}{\begin{eqnarray}}
\newcommand{\eea}{\end{eqnarray}} \newcommand{\ee}{\end{equation}}
 
 \def\ba{\begin{eqnarray}}
\def\ea{\end{eqnarray}}

\def\L{{\Lambda}}

\def\log{\,{\rm log}\,}

\def\p{\partial}

\def\a{\alpha}

\def\e{\epsilon}

\def\[{{\bf [}}
\def\]{{\bf ]}}

\begin{document}
\starttext \baselineskip=18pt \setcounter{footnote}{0}
\newtheorem{theorem}{Theorem}
\newtheorem{lemma}{Lemma}
\newtheorem{definition}{Definition}
\begin{center} {\Large \bf ON THE SCALING LIMIT OF A SINGULAR INTEGRAL OPERATOR}
\footnote{Research supported in part by National Science Foundation
grants  DMS-02-45371 and DMS-04-05519}
\\
\bigskip

{\large I.M. Krichever and D.H. Phong} \\

\bigskip

Department of Mathematics\\
Columbia University, New York, NY 10027

\end{center}

\baselineskip=15pt \setcounter{equation}{0} \setcounter{footnote}{0}

\section{Introduction}
\setcounter{equation}{0}

Recently an alternative approach to Birkhoff's theory of difference equations
\cite{B} has been proposed in \cite{K}. This approach leads naturally to local
monodromies of difference equations, which should converge in principle to
monodromy matrices of differential equations, thus providing a missing link
in the theory of isomonodromic transformations of systems of linear difference
equations (see e.g. \cite{Bo, HI, TW} and references therein).

\medskip
The key to the convergence process in \cite{K} is the scaling
limit of a certain singular integral operator $I$, arising from
a Riemann-Hilbert problem.
The operator $I$ acts on functions $\phi$ defined on the vertical line
with fixed abscisse at $a$, and its kernel $k(z,\xi)$ is given explicitly by
\be
k(z,\xi)={e^{\pi i(z-a)}+e^{-\pi i(z-a)}\over
(e^{i\pi(\xi-a)}+e^{-\pi i(\xi-a)})
(e^{i\pi(\xi-z)}+e^{-\pi i(\xi-z)})}
\ee
If we set $\pi z=\pi a+iy$, $\pi \xi=\pi a+i\eta$, $y,\eta\in{\bf R}$,
and view $\phi$ as a function of $\eta$, we can define the following re-scaled
versions $I_\lambda$ of $I$,
\be
I_\lambda(\phi)(y)=P.V.\int_{-\infty}^\infty\ {1\over e^{\lambda(y-\eta)}-e^{-\lambda(y-\eta)}}
{e^{-\lambda y}+e^{\lambda y}\over e^{-\lambda\eta}+e^{\lambda\eta}}\ \phi(\eta)\,d\eta,
\ee
where $P.V.$ denotes principal values. As noted in \cite{K}, an essential property of
the operators $I_\lambda$ is their formal limit,
\be
\label{limit0}
I_\lambda(\phi)(y)\to\int_0^y\phi(\xi)d\xi,
\qquad\lambda\to+\infty.
\ee

The purpose of the present paper is to provide a detailed study of the boundedness
properties of the operators $I$ and $I_\lambda$ in suitable spaces of Schauder type,
and to establish a precise version of the formal limit (\ref{limit0}).
Near the diagonal, the singularities of the kernels of $I_\lambda$ are the same as for
the Hilbert transform, and the techniques for handling the local behavior
of such kernels
are well-known. The main novel feature in our case is rather their global behavior near
$\infty$. This global behavior prevents their boundedness on scale-invariant spaces,
and accounts for the existence of non-trivial limits such as (\ref{limit0}).

\section{Schauder estimates with exponential
growth}
\setcounter{equation}{0}

We introduce the following norms of Schauder type for functions on ${\bf R}$.
Fix $\kappa\in{\bf R}$, $m\in {\bf Z}$, $0<\a<1$, and let $\Lambda_{(m,\kappa)}^\a$
be the space of functions $\phi$ on ${\bf R}$ satisfying the conditions
\bea
\label{norm}
|\phi(x)|&\leq& C\,(1+|x|)^m\,e^{\kappa|x|},\nonumber\\
|\phi(x)-\phi(y)|
&\leq&
C\,|x-y|^\a\,\big\{(1+|x|)^m\,e^{\kappa|x|}+(1+|y|)^me^{\kappa|y|}\big\},
\eea
for all $x,y\in{\bf R}$. We define $||\phi||_{\Lambda_{(m,\kappa)}^\a}$
to be the infimum of the constants $C$ for which these inequalities hold.
We also require the space $\Lambda_{(log,\kappa)}^\a$ and the corresponding norm
$||\phi||_{\Lambda_{(log,\kappa)}^\a}$ defined by the conditions
\bea
\label{norm1}
|\phi(x)|&\leq& C\,\log\,(1+|x|)\,e^{\kappa|x|},\nonumber\\
|\phi(x)-\phi(y)|
&\leq&
C\,|x-y|^\a\,\big\{\log\,(1+|x|)\,e^{\kappa|x|}+\log\,(1+|y|)e^{\kappa|y|}\big\},
\eea

\medskip
The singular integral operator $I$ can be expressed as
\be
I(\phi)(y)=(e^{-y}+e^y)H({1\over e^{-(\cdot)}+e^{(\cdot)}}\phi(\cdot))
\ee
where $H$ is the following exponentially decaying version of the
classical Hilbert transform,
\be
\label{hilbert}
(H\psi)(y)=P.V.\int_{-\infty}^\infty {1\over e^{y-\eta}-e^{-(y-\eta)}}\,\psi(\eta)\,d\eta
\equiv
{\rm lim}_\e\to 0
\int_{|y-\eta|>\e} {1\over e^{y-\eta}-e^{-(y-\eta)}}\,\psi(\eta)\,d\eta.
\ee
Set
\be
\label{hilbertkernel}
K(z)={1\over e^z-e^{-z}}.
\ee
Then the kernel $K(z)$ is $C^\infty({\bf R}\setminus 0)$, odd, and satisfies
\be
|K(z)|\leq C\,\cases{|z|^{-1}, &if $|z|\leq 1$;\cr
e^{-|z|}, &if $|z|>1$.\cr}
\qquad
|\p_zK(z)|\leq C\,\cases{|z|^{-2}, &if $|z|\leq 1$;\cr
e^{-|z|}, &if $|z|>1$.\cr}
\ee
In particular, these are better estimates than for the standard Hilbert transform kernel
$K_0(z)=z^{-1}$, and it follows at once that the operator $H$ is bounded on the standard
Schauder spaces (see e.g. \cite{Stein, LU}). To obtain estimates for the operator $I$,
we need the boundedness of $H$ on the above spaces $\Lambda_{(m,\kappa)}^\a$,
and this is provided by the following theorem:

\begin{theorem}
\label{Hboundtheorem}

Fix $0<\a<1$, $m\in {\bf Z}$.
The operator $H$ is bounded on the following Schauder spaces,
\be
\label{Hbound}
||H\psi||_{\Lambda_{(m,\kappa)}^\a}
\leq
\ C_{m,\a,k}\, ||\psi||_{\Lambda_{(m,\kappa)}^{\a}},
\qquad -1<\kappa<1.
\ee
For $\kappa=-1$, we have the following bounds, for $m\in{\bf Z}$, $m\geq-1$,
\bea
\label{Hbound1}
||H\psi||_{\Lambda_{(m+1,-1)}^\a}
&\leq&
\ C_{m,\a}\, ||\psi||_{\Lambda_{(m,-1)}^{\a}},
\qquad m\geq 0,
\nonumber\\
||H\psi||_{\Lambda_{(log,-1)}^\a}
&\leq&
\ C_{\a}\, ||\psi||_{\Lambda_{(m,-1)}^{\a}},
\qquad m=-1.
\eea

\end{theorem}

\bigskip

\noindent
{\it Proof.} The method of proof is the standard method for Schauder estimates
for singular integral operators. The only new feature here is the control of $H\psi(x)$
for $x$ large. In view of the fact that $K(z)$ is odd
and exponentially decreasing, we can write
\be
H\psi(x)
=
\int_{-\infty}^\infty K(x-y)\,(\psi(y)-\psi(x)) dy,
\ee
where the integrals on the right hand side are now convergent
for $\psi\in\Lambda_{(m,\kappa)}^\a$ with $0<\a<1$, $\kappa<1$.
In particular,
\bea
|H\psi(x)|
&\leq&
\int_{|x-y|<1}
|x-y|^{-1+\a}
((1+|x|)^m e^{\kappa|x|}+(1+|y|)^m e^{\kappa|y|})\,dy
\nonumber\\
&&
+
\int_{|x-y|\geq 1}
e^{-|x-y|}((1+|x|)^m e^{\kappa|x|}+(1+|y|)^m e^{\kappa|y|})\,dy.
\eea
These are clearly bounded for $|x|$ bounded, so we may assume that $|x|\geq 3$. In this case,
${1\over 2|x|}\leq |x|-1\leq |y|\leq |x|+1\leq 2|x|$ in the integral over the region $|x-y|<1$, and
$(1+|y|)^me^{\kappa|y|}\leq C_\kappa\,(1+|x|)^{m}e^{\kappa|x|}$.
Thus the first integral is bounded by $C\,(1+|x|)^{m}e^{\kappa|x|}$.
The same upper bound for the second integral follows from the following lemma:

\begin{lemma}
\label{integralestimates1}
For any $-1<\kappa<1$, and any $m\in{\bf Z}$, we have for all $|x|>3$
\be
\label{integrals1}
\int_{\bf R} e^{-|z|}(1+|x-z|)^me^{\kappa|x-z|}\,dz
\leq
C_{m,\kappa}(1+|x|)^m e^{\kappa |x|}.
\ee
For $\kappa=-1$, we have for $m\in {\bf Z}$, $m\geq -1$,
\be
\label{integrals2}
\int_{\bf R} e^{-|z|}(1+|x-z|)^me^{-|x-z|}\,dz
\leq
C_{m}\cases{e^{-|x|}(1+|x|)^{m+1}, &if $m\geq0$\cr
e^{-|x|}\log\,(1+|x|), &if $m=-1$.\cr}
\ee
\end{lemma}

\noindent
{\it Proof of Lemma \ref{integralestimates1}.} We consider separately the cases
of $0\leq \kappa<1$, $-1<\kappa<0$, and $\kappa=-1$. When $0\leq\kappa<1$, we write
$e^{\kappa|x-z|}\leq e^{\kappa|x|}\,e^{\kappa|z|}$, and hence
the integral on the left hand side of the above inequality can be bounded by
\be
e^{\kappa|x|}\int_{|x-z|\geq {1\over 2}|x|}e^{-(1-\kappa)|z|}(1+|x-z|)^m dz
+
e^{\kappa|x|}\int_{|x-z|< {1\over 2}|x|}e^{-(1-\kappa)|z|}(1+|x-z|)^m dz.
\ee
In the first integral we can write
\be
(1+|x-z|)^m\leq C_m (1+|x|)^m(1+|z|)^m.
\ee
This is certainly true with $C_m=1$ if $m\geq 0$. If $m<0$, then we use the condition
$|x-z|\geq {1\over 2}|x|$ to write
$(1+|x-z|)^m\leq 2^{-m}(1+|x|)^m$, and the inequality still holds.
Since $\kappa<1$, the desired bound follows for the first integral.
Next, in the second integral, we have ${1\over 2}|x|<|z|<{3\over 2}|x|$, and
we can write
\bea
\int_{|x-z|< {1\over 2}|x|}e^{-(1-\kappa)|z|}(1+|x-z|)^m dz
&\leq&
e^{-{1-\kappa\over 2}|x|}
\int_{|z|\leq{3\over 2}|x|}
(1+|x|)^{|m|}(1+|z|)^{|m|}\,dz
\nonumber\\
&\leq& C_N\,(1+|x|)^{-N},
\eea
for arbitrary $N$. This proves the lemma when $0\leq\kappa<1$.
When $-1<\kappa<0$, we write instead
\be
e^{-|z|}e^{\kappa|x-z|}
=
e^{-(1+\kappa)|z|}e^{\kappa(|z|+|x-z|)}
\leq
e^{-(1+\kappa)|z|}e^{\kappa|x|}
\ee
and bound the integral on the left hand side of the lemma by
\be
e^{\kappa|x|}
\int_{|x-z|\geq{1\over 2}||x|}e^{-(1+\kappa)|z|}(1+|x-z|)^mdz
+
e^{\kappa|x|}
\int_{|x-z|<{1\over 2}||x|}e^{-(1+\kappa)|z|}(1+|x-z|)^mdz.
\ee
The bounds for these integrals are now the same as in the previous case.
This establishes the estimate (\ref{integrals1}).
Finally, consider the case $\kappa=-1$.
In the region of integration $|x-z|>4|x|$, we have
the integrand can be crudely bounded by
$e^{-2|x|}(1+|x-z|)^m e^{-{1\over 2}|x-z|}$,
and hence the contribution of this region is $O(e^{-2|x|})$,
which is better than we actually need.
Thus it suffices to consider the region $|x-z|\leq 4|x|$.
We write then
\be
\int_{|x-z|<4|x|}e^{-|z|}(1+|x-z|)^m e^{-|x-z|}dz
\leq
e^{-|x|}\int_{|x-z|<4|x|}(1+|x-z|)^m
\ee
from which the desired estimate follows at once.
The proof of the lemma is complete. Q.E.D.

\medskip
We return to the proof of the theorem. Let $x,x'\in{\bf R}$
and set $\delta=|x-x'|$. The next step is to estimate
$H\psi(x)-H\psi(x')$, which can be expressed as
\bea
\label{breakup1}
&&
\int_{|y-x|<3\delta}K(x-y)(\psi(y)-\psi(x))dy
-
\int_{|y-x'|<3\delta}K(x'-y)(\psi(y)-\psi(x'))dy
\\
&&
+
\int_{|y-x|\geq 3\delta}K(x-y)(\psi(y)-\psi(x))dy
-
\int_{|y-x'|\geq 3\delta}K(x'-y)(\psi(y)-\psi(x'))dy.
\nonumber
\eea
The first two integrals can be estimated as in the bounds for $|H\psi(x)|$. For example,
\bea
|\int_{|y-x|<3\delta}K(x-y)(\psi(y)-\psi(x))dy|
&\leq&
||\psi||_{\Lambda_{(m,\kappa)}^\a}
\int_{|x-y|<3\delta}
|x-y|^{-1+\a}
\big\{(1+|x|)^me^{\kappa|x|}
\nonumber\\
&&
\qquad\qquad\qquad
+
(1+|y|)^me^{\kappa|y|}\big\}
\nonumber\\
&\leq&
C\,||\psi||_{\Lambda_{(m,\kappa)}^\a}\,\delta^\a(1+|x|)^me^{\kappa|x|}
\eea
since $(1+|y|)^me^{\kappa|y|}\leq C\,(1+|x|)^me^{\kappa|x|}$
for $|x|\geq 3$ and $\delta<<1$. To estimate the remaining two
integrals, write
\bea
\int_{|y-x'|\geq 3\delta}K(x'-y)(\psi(y)-\psi(x'))dy
&=&
\int_{|y-x'|\geq 3\delta}K(x'-y)(\psi(y)-\psi(x))dy
\\
&=&
\int_{|y-x|\geq 3\delta}
+
\int_{|y-x'|\geq3\delta,
|y-x|<3\delta}
-
\int_{|y-x'|<3\delta,|y-x|>3\delta}
\nonumber
\eea
The last two integrals on the right hand side satisfy the desired bounds,
because in their ranges of integration, we have $|y-x|\sim |y-x'|\sim\delta$,
and the same arguments above apply. The remaining integral can be combined with
the third integral in (\ref{breakup1}) to give
\bea
\int_{|y-x|>3\delta}
(K(x-y)-K(x'-y))(\psi(y)-\psi(x))dy.
\eea
Since we have
\be
|K(x-y)-K(x'-y)|\leq |x-y|\cdot |\p_xK(z)|
\ee
for some $z$ in the segment between $x-y$ and $x'-y$,
and hence $|z|\sim|x-y|$ when $|y-x|>3|x-x'|$, we can write,
in view of the bounds for the $|\p_zK(z)|$,
\bea
&&
|\int_{|y-x|>3\delta}
(K(x-y)-K(x'-y))(\psi(y)-\psi(x))dy|
\nonumber\\
&&
\quad
\leq
\delta\,||\psi||_{\Lambda_{m,\kappa}^\a}
\int_{3\delta<|x-y|<1}|x-y|^{-2+\a}
\big\{(1+|x|)^me^{\kappa|x|}
+
(1+|y|)^me^{\kappa|y|}\big\}\,dy
\nonumber\\
&&
\quad\quad
+
\delta\,
||\psi||_{\Lambda_{m,\kappa}^\a}
\int_{|x-y|>1}e^{-|x-y|}
\big\{(1+|x|)^me^{\kappa|x|}
+
(1+|y|)^me^{\kappa|y|}\big\}\,dy.
\eea
The first integral on the right hand side is bounded by
\bea
&&
\delta\,||\psi||_{\Lambda_{m,\kappa}^\a}
\int_{3\delta<|x-y|<1}|x-y|^{-2+\a}
\big\{(1+|x|)^me^{\kappa|x|}
+
(1+|y|)^me^{\kappa|y|}\big\}\,dy\nonumber\\
&&
\quad
\leq
\delta\,||\psi||_{\Lambda_{m,\kappa}^\a}\,
(1+|x|)^m e^{\kappa|x|}
\int_{3\delta<|x-y|<1}|x-y|^{-2+\a}
\leq C\,\delta^\a \,(1+|x|)^m e^{\kappa|x|}.
\eea
Applying Lemma (\ref{integralestimates1}), we obtain similar bounds for
the second integral. Altogether, we have shown that
\be
|H\psi(x)-H\psi(x')|
\leq C\,||\psi||_{\Lambda_{(m,\kappa)}^\a}
|x-x'|^\a (1+|x|)^me^{\kappa|x|}
\ee
for $|x-x'|$ small, and the theorem is proved when $-1<\kappa<1$.
The case $\kappa=-1$ is established exactly in the same way,
using the corresponding estimates in Lemma \ref{integralestimates1}
for $\kappa=-1$ and $m\geq -1$. The proof of Theorem \ref{Hboundtheorem}
is complete. Q.E.D.

\bigskip

\begin{theorem}
\label{bound}
For $0<\kappa<2$, $m\in{\bf Z}$, and
$0<\a<1$, the operator $I$ is bounded on the following Schauder spaces,
\bea
||I(\phi)||_{\Lambda_{(m,\kappa)}^\a}
\leq\ C_{m,\kappa,\a}\, ||\phi||_{\Lambda_{(m,\kappa)}^{\a}},
\qquad
0<\kappa<2.
\eea
For $\kappa=0$, $m\in{\bf Z}$, $m\geq -1$,
the operator $I$ satisfies the following bounds,
\bea
||I(\phi)||_{\Lambda_{(m+1,\kappa)}^\a}
&\leq&\ C_{m,\a}\, ||\phi||_{\Lambda_{(m,\kappa)}^{\a}},
\qquad
m\geq 0
\nonumber\\
||I(\phi)||_{\Lambda_{(log,\kappa)}^\a}
&\leq&\ C_{\a}\, ||\phi||_{\Lambda_{(m,\kappa)}^{\a}},
\qquad
m=-1.
\eea

\end{theorem}

\bigskip
\noindent
{\it Proof.} This is an easy consequence of Theorem 1,
the fact that the
map $\phi\to\psi(y)={1\over e^y+e^{-y}}\phi(y)$ is a one-to-one
and onto map $M$ from $\Lambda_{(m,\kappa)}^\a\to \Lambda_{(m,\kappa-1)}^\a$,
with equivalent norms
\be
||\psi||_{\Lambda_{(m,\kappa-1)}^\a}
\sim
||\phi||_{\Lambda_{(m,\kappa)}^\a}.
\ee
and the relation $I(\phi)=M^{-1}HM\phi$. Q.E.D.

\bigskip
We observe that these bounds always require some space which is not scale-invariant.
Thus bounds for $I_\lambda$ cannot be obtained by scaling the bounds for $I$, and this
explains partly the possibility of the scaling limits discussed in the next section.

\section{The scaling limit of $I_\lambda$}
\setcounter{equation}{0}

We come now to the operators $I_\lambda$. The estimates for $I$
in the previous section show that $I_\lambda$ cannot be treated by
simple scaling arguments from $I$. Instead, we shall study the bounds
and limits for $I_\lambda$ as $\lambda\to+\infty$ directly. It is simplest
to carry this out for functions $\phi$
satisfying conditions of the form,
\bea
|\p^l\phi(x)|
\leq&
C_k\,(1+|x|)^m,
\qquad 0\leq l\leq k,
\eea
for fixed $m\in {\bf N}$, $k\in {\bf N}$, and norms $||\phi||_{\Lambda_{(m)}^k}$
defined to be the best constant $C_k$ for which the above condition holds.
The following theorem describes the limit of $I_\lambda$ in these spaces, although
it should be clear from the proof and from the previous section that
other more precise versions can be formulated as well:

\begin{theorem}
\label{limit}
Fix $m\in {\bf N}$. Then we have the following bounds, uniform in $\lambda$ and in
$\phi\in \Lambda_{(m)}^{1}$,
\be
\label{mainestimate}
||I_\lambda(\phi)(y)
-
\int_0^y\phi(\xi)d\xi||_{\Lambda_{(m+1)}^0}
\leq
C_{m}\, \lambda^{-{1\over 2}}||\phi||_{\Lambda_{(m)}^{1}}.
\ee

\end{theorem}

\bigskip
\noindent
{\it Proof.} Formally, if we write
\be
I_\lambda\phi(y)=\int K_\lambda(y,\eta)\phi(\eta)\,d\eta
\ee
with
\be
K_\lambda(y,\eta)
=
{e^{-\lambda y}+e^{\lambda y}\over e^{-\lambda\eta}+e^{\lambda\eta}}
{1\over e^{\lambda(y-\eta)}-e^{-\lambda(y-\eta)}}
\ee
then for, say, $y>0$, we have the pointwise limit
\be
K_\lambda(y,\eta)
\to \cases{1,\ \ \ &if $0<\eta<y$\cr
0,\ \ \ &if $\eta<0$ or $\eta>y$.\cr}
\ee
Thus, formally, the left hand side of the expression in the theorem tends to $0$ as
$\lambda\to+\infty$. However, none of the integrals involved is uniformly nor absolutely
convergent, and we have to proceed with care.
Fix $y>0$ (the case of $y<0$ being similar). The key to the estimates
is the following break-up of the principal value integral defining $I_\lambda(\phi)$,
\bea
\label{breakup}
(I_\lambda(\phi))(y)
&=&
\int_0^y
{e^{\lambda y}+e^{-\lambda y}\over e^{\lambda t}-e^{-\lambda t}}\,(\,\psi(y-t)-\psi(y+t)\,)dt
+
\int_{|t|>y}{e^{\lambda y}+e^{-\lambda y}\over e^{\lambda t}-e^{-\lambda t}}\psi(y-t)\,dt
\nonumber\\
&\equiv&
(A)+(B)
\eea
with
\be
\psi(\eta)={1\over e^{\lambda\eta}+e^{-\lambda\eta}}\phi(\eta).
\ee
To estimate $(A)$, we apply Taylor's formula
\bea
\psi(y-t)-\psi(y+t)
=t\int_{-1}^1 \psi'(y-\rho t)d\rho
\eea
which gives in this particular case,
\bea
\psi(y-t)-\psi(y+t)
&=&
t\int_{-1}^1\,d\rho\,(\ {1\over e^{\lambda(y-\rho t)}+e^{-\lambda(\eta-\rho t)}}\phi'(y-\rho t)
\nonumber\\
&&
\qquad\qquad
-\lambda {e^{\lambda(y-\rho t)}-e^{-\lambda(y-\rho t)}\over
(e^{\lambda(y-\rho t)}+e^{-\lambda(y-\rho t)})^2}\phi(y-\rho t)\ )
\eea
Thus $(A)$ can be rewritten as
\bea
(A)&=&
{1\over\lambda}
\int_0^y dt\int_{-1}^1 d\rho
\chi_\lambda(\rho,t)
\phi'(y-\rho t)
\nonumber\\
&&
-\int_0^y dt\int_{-1}^1 d\rho
\chi_\lambda(\rho,t)\,
{e^{\lambda(y-\rho t)}-e^{-\lambda(y-\rho t)}\over
e^{\lambda(y-\rho t)}+e^{-\lambda(y-\rho t)}}\phi(y-\rho t)
\nonumber\\
&\equiv& A_1+A_0
\eea
where the function $\chi_\lambda(\rho,t)$ is defined by
\bea
\label{chi}
\chi_\lambda(\rho,t)
=
{\lambda t\over e^{\lambda t}-e^{-\lambda t}}\,{e^{\lambda y}+e^{-\lambda y}\over
e^{\lambda(y-\rho t)}+e^{-\lambda (y-\rho t)}}.
\eea
The following sharp estimates for $\chi_\lambda(\rho,t)$ play an essential
role in the sequel:

\begin{lemma}
\label{lemmadirac}

For all $0<t<y$,
the functions $\chi_\lambda(\rho,t)$ satisfy the following properties
\bea
\label{dirac}
&&
{\rm (a)}\qquad
{1\over 2}\,
{\lambda t\over 1-e^{-2\lambda t}}
\,e^{-\lambda t(1-\rho)}
<\chi_\lambda(\rho,t)\leq \,2\,{\lambda t\over 1-e^{-2\lambda t}}
\,e^{-\lambda t(1-\rho)}, \qquad |\rho|<1,
\ 0<t<y \nonumber\\
&&
{\rm (b)}\qquad {1\over 2}\leq \int_{-1}^1\chi_\lambda(\rho,t)d\rho \leq 2.
\eea
\end{lemma}

\noindent
{\it Proof.} In the region $0<t<y$, we have
$y-\rho t >0$ for all $|\rho| <1$,
and thus
\bea
{1\over 2}\,
e^{\lambda \rho t}
\leq {e^{\lambda y}+e^{-\lambda y}\over e^{\lambda(y-\rho t)}+e^{-\lambda(y-\rho t)}}
\leq
2\,e^{\lambda\rho t}.
\eea
The upper bound implies (a), while the lower bound implies (b), when combined
with the following explicit formula
\be
\int_{-1}^1e^{\lambda\rho t}\,d\rho={1\over\lambda t}(\ e^{\lambda t}-e^{-\lambda t}\ ).
\ee
The proof of Lemma \ref{lemmadirac} is complete.

\medskip
We can now show that $A_1\to 0$ with a precise rate:

\begin{lemma}
\label{lemmaphiderivative}
The term involving $\phi'$ above tends to $0$ at the following rate,
\bea
|A_1|\leq C_m\,{1\over\lambda}\,||\phi||_{\Lambda_{(m,0)}^1}(1+y)^{m+1}.
\eea
\end{lemma}

\noindent
{\it Proof}. It suffices to write
\be
|A_1|
\leq
\,{1\over\lambda}\,
||\phi||_{\Lambda_{(m)}^1}\,\int_0^y dt\,(1+|y|^m+|t|^m)\,
\int_{-1}^1d\rho\ \chi_\lambda(\rho,t)
\ee
and the desired estimate follows from the statement (b)
of Lemma \ref{lemmadirac}. Q.E.D.

\bigskip
The estimates in Lemma \ref{lemmadirac} show
that $\chi_\lambda(\rho,t)$ provide an approximation of the Dirac measure
concentrated at $\rho=1$,
\be
{1\over\int_{-1}^1\chi_\lambda(\mu,t)d\mu}\,\chi_\lambda(\rho,t)\to \delta(\mu-1)
\ee
A precise version of this statement with sharp estimates
is given in the next lemma. Set
\bea
\label{chibreakup}
&&
\int_{-1}^1 \,d\rho \,\chi_\lambda(\rho,t)\,
{e^{-\lambda(y-\rho t)}-e^{\lambda(y-\rho t)}
\over
e^{-\lambda(y-\rho t)}+e^{\lambda(y-\rho t)}}\phi(y-\rho t)
-\phi(y-t)
\\
&&
\quad
=
\int_{-1}^1 \,d\rho \,\chi_\lambda(\rho,t)\,
({e^{-\lambda(y-\rho t)}-e^{\lambda(y-\rho t)}
\over
e^{-\lambda(y-\rho t)}+e^{\lambda(y-\rho t)}}-1)
\,\phi(y-\rho t)
\nonumber\\
&&
\quad\quad
+
\int_{-1}^1 \,d\rho \,\chi_\lambda(\rho,t)\,(\phi(y-\rho t)-\phi(y-t))
+
(\int_{-1}^1 \,d\rho \,\chi_\lambda(\rho,t)-1)\phi(y-t).
\nonumber
\eea

\begin{lemma}
\label{lemmadiracestimates}
For all $0<t<y$, and any $\delta>0$ and small,
we have the following estimates, with absolute constants,
\bea
\label{diracestimates}
&&
{\rm (a)}\quad
\bigg|\int_{-1}^1 \,d\rho \,\chi_\lambda(\rho,t)-1\bigg|\cdot
|\phi(y-t)|
\leq
\,||\phi||_{\Lambda_{(m)}^0}(1+y)^m(e^{-\lambda y}+e^{-2\lambda (y-t)})
\\
&&
{\rm (b)}\quad
\bigg|\int_{-1}^1 \,d\rho \,\chi_\lambda(\rho,t)\,(\phi(y-\rho t)-\phi(y-t))\bigg|
\leq \,||\phi||_{\Lambda_{(m)}^1}(1+y)^m(\ \delta t+e^{-\lambda\delta t}\ )
\\
&&
{\rm (c)}\quad
\int_{-1}^1\,d\rho\,\chi_\lambda(\rho,t)\bigg|
{e^{\lambda(y-\rho t)}-e^{-\lambda(y-\rho t)}
\over
e^{\lambda(y-\rho t)}+e^{-\lambda(y-\rho t)}}-1
\bigg|\,|\phi(y-\rho t)|
\leq \, C\,e^{-2\lambda(y-t)}\,||\phi||_{\Lambda_{(m)}^0}(1+y)^m.
\nonumber
\\
\eea
\end{lemma}

\noindent
{\it Proof.} To prove (a), we write
\bea
\bigg|{e^{-\lambda y}+e^{\lambda y}
\over
e^{-\lambda(y-\rho t)}-e^{\lambda(y-\rho t)}}
-
e^{\lambda\rho t}
\bigg|
&=&
e^{\lambda\rho t}
\bigg|
{e^{-2\lambda y}-e^{-2\lambda(y-\rho t)}
\over
1+e^{-2\lambda(y-\rho t)}}
\bigg|
\\
&\leq&
e^{\lambda \rho t}\big(e^{-2\lambda y}+e^{-2\lambda(y-\rho t)})
\leq
e^{\lambda \rho t}\big(e^{-2\lambda y}+e^{-2\lambda(y-t)}).
\nonumber
\eea
In particular,
\be
\bigg|
\int_{-1}^1 \,d\rho \,\chi_\lambda(\rho,t)
-
{\lambda t\over e^{\lambda t}-e^{-\lambda t}}
\int_{-1}^1\,d\rho\,e^{\lambda\rho t}
\bigg|
\leq
{\lambda t\over e^{\lambda t}-e^{-\lambda t}}
\int_{-1}^1\,d\rho\,e^{\lambda\rho t}
(e^{-2\lambda y}+e^{-2\lambda (y-t)}).
\ee
Since $\int_{|\rho|<1}\,d\rho\,e^{\lambda\rho t}
=(\lambda t)^{-1}(e^{\lambda t}-e^{-\lambda t})$,
the statement (a) follows.

\medskip

To establish the statement (c), we begin by noting that
\be
e^{-2\lambda(y-\rho t)}
\leq 1-{e^{\lambda(y-\rho t)}-e^{-\lambda(y-\rho t)}
\over
e^{\lambda(y-\rho t)}+e^{-\lambda(y-\rho t)}}
=
2 \,{e^{-2\lambda(y-\rho t)}
\over
1+e^{-2\lambda(y-\rho t)}}
\leq 2\,e^{-2\lambda(y-\rho t)}.
\ee
Using the estimate for $\chi_\lambda$ in Lemma \ref{lemmadirac}
and carrying out explicitly the integral in $\rho$ gives
\bea
\int_{-1}^1\,d\rho\,\chi_\lambda(\rho,t)\bigg|
{e^{\lambda(y-\rho t)}-e^{-\lambda(y-\rho t)}
\over
e^{\lambda(y-\rho t)}-e^{-\lambda(y-\rho t)}}-1
\bigg|
&\leq&
{\lambda t\over 1-e^{-2\lambda t}}e^{-\lambda(t+2y)}\int_{-1}^1\,d\rho
e^{3\lambda \rho t}
\nonumber\\
&=&
{1\over 3}e^{-2\lambda(y-t)}{1-e^{-6\lambda t}\over 1-e^{-2\lambda t}}
\leq C\, e^{-2\lambda(y-t)},
\nonumber
\eea
which implies immediately (c).

\medskip

To establish (b), let $\delta>0$ be any number sufficiently small
and to be chosen suitably later. Write
\bea
\int_{-1}^1 \,d\rho \,\chi_\lambda(\rho,t)\,(\phi(y-\rho t)-\phi(y-t))
=
\int_{-1}^{1-\delta}+\int_{1-\delta}^1
\equiv I_\delta +II_{\delta}.
\eea
The second term on the right hand side can be estimated by,
\be
|II_\delta|
\leq
\delta\,t\,||\phi||_{\Lambda_{(m)}^1}(1+y)^m\int_{1-\delta}^1\,d\rho\,\chi_\lambda(\rho,t)
\leq
\delta\,t\,||\phi||_{\Lambda_{(m)}^1}(1+y)^m,
\ee
while the first term can be estimated using Lemma \ref{lemmadirac},
\bea
|I_\delta|
&\leq&
2\,||\phi||_{\Lambda_{(m)}^1}(1+y)^m\,
{\lambda t\over 1-e^{-2\lambda t}}
\int_{-1}^\delta\,d\rho\,\,e^{-\lambda t(1-\rho)}
\nonumber\\
&=&
2\,||\phi||_{\Lambda_{(m)}^1}(1+y)^m\,e^{-\lambda\delta t}
{1-e^{-\lambda t(2-\delta)}\over 1-e^{-2\lambda t}}
\leq
C\,{\rm sup}_{[0,2y]}|\phi|\,e^{-\lambda\delta t}.
\eea
The proof of Lemma \ref{lemmadiracestimates} is complete. Q.E.D.

\bigskip
We can now carry out the integral in $t$. The precise estimates
are given in the next lemma:

\begin{lemma}
\label{lemmatintegral0}
For any $0<y$, we have the following estimates,
\bea
&&
\bigg|
\int_0^y dt\int_{-1}^1 d\rho {\lambda \,t\over e^{\lambda t}-e^{-\lambda t}}
(e^{\lambda y}+e^{-\lambda y}){e^{\lambda(y-\rho t)}-e^{-\lambda(y-\rho t)}\over
(e^{\lambda(y-\rho t)}+e^{-\lambda(y-\rho t)})^2}\phi(y-\rho t)
-
\int_0^y \,dt\,\phi(y-t)
\bigg|
\nonumber\\
&&
\qquad\qquad\leq
C\,
||\phi||_{\Lambda_{(m)}^1}(1+y)^m\,(\ {y\over 1+\lambda y}+{y\over\lambda^{1\over 2}}\ ).
\eea
with a constant $C$ independent of $y$ and of $\lambda$.

\end{lemma}

\noindent
{\it Proof.} In view of the defining formula (\ref{chi}) for
the function $\chi_\lambda(\rho,t)$ and the break up (\ref{chibreakup}),
the left hand side of the desired inequality
is bounded by the integral in $t$ of the three inequalities
in Lemma \ref{lemmadiracestimates}. This gives the following upper bound,
\bea
||\phi||_{\Lambda_{(m)}^1}(1+y)^m\ \int_0^y dt\,
(e^{-\lambda y}+2 e^{-2\lambda(y-t)}
+\delta t+e^{-\delta \lambda t}).
\eea
The integral can be evaluated explicitly, and we find
\be
ye^{-\lambda y}+{1\over\lambda}(1-e^{-2\lambda y})
+
{1\over 2}\delta y^2+{1\over\delta\lambda}(1-e^{-\delta\lambda}).
\ee
We consider the sum of the first two terms: when $\lambda y<1$,
it is bounded by $C\,y$, where $C$ is an absolute constant.
When $\lambda y\geq 1$, it is bounded by $C\,\lambda^{-1}$. Thus we have
\be
ye^{-\lambda y}+{1\over\lambda}(1-e^{-2\lambda y})\sim {y\over 1+\lambda y}.
\ee
Next, we consider the optimal choice of $\delta$ so as to minimize the
size of the sum of the remaining two terms in the above integral. We note
that we may assume that $\delta\lambda >1$, since otherwise
the term $(\delta\lambda)^{-1}(1-e^{-\delta\lambda})$ is of size $1$,
and we do not even get convergence to $0$. Thus we should take $\delta\lambda>1$,
in which case the sum of the two remaining terms is of size
\be
\delta y^2+{1\over\delta\lambda}
\ee
which attains its lowest size $y\lambda^{-{1\over 2}}$ if we set $\delta=y^{-1}\lambda^{-{1\over 2}}$.
This gives the estimate stated in the lemma. Q.E.D.

\bigskip
We return now to the estimate of the contribution to $I_\lambda(\phi)(y)$
of the integral in $t$ from the region $|t|>y$.

\begin{lemma}
\label{lemmatinfinity}
For any $0<y$, we have the following estimate
\bea
\bigg|\int_{|t|>y}{1\over e^{\lambda t}-e^{-\lambda t}}
{e^{\lambda y}+e^{-\lambda y}\over
e^{\lambda(y-t)}+e^{-\lambda(y-t)}}
\phi(y-t)dt
\bigg|
&\leq&
C_m\,{1\over\lambda}||\phi||_{\Lambda_{(m)}^0}(1+\log (1+{1\over \lambda y})
\nonumber
\eea
\end{lemma}

\noindent
{\it Proof.} Consider first the contribution from the region $t>y$. In this region,
we have
\be
{1\over 2}e^{\lambda(2y-t)}
\leq
{e^{\lambda y}+e^{-\lambda y}\over e^{\lambda(y- t)}+e^{-\lambda(y- t)}}
\leq
2\,{e^{\lambda y}\over e^{\lambda(t-y)}}=
2\,e^{\lambda(2y-t)}
\ee
Thus the contribution from the region $t>y$
to the integral
on the left hand side of the desired inequality
can be bounded by
\bea
\int_{t>y}
{e^{-2\lambda(t-y)}\over 1-e^{-2\lambda t}}|\phi(y-t)|\,dt
&=&
\int_0^\infty
{e^{-2\lambda s}\over 1-e^{-2\lambda(s+y)}}|\phi(-s)|\,ds
\nonumber\\
&\leq &
||\phi||_{\L_{(m)}^0}
\int_0^\infty
{e^{-2\lambda s}\over 1-e^{-2\lambda(s+y)}}(1+|s|^m)\,ds.
\eea
We claim that for all $m\in {\bf N}$, we have
\be
\int_0^\infty
{e^{-2\lambda s}\over 1-e^{-2\lambda(s+y)}}(1+|s|^m)\,ds
\leq
C_m\,
{1\over\lambda}(1+\log (1+{1\over \lambda y}).
\ee
In fact, setting $\mu=e^{-2\lambda y}$ and making the change of variables
$s\to u$, $e^{-2\lambda u}=s$, this integral can be rewritten as
\be
{1\over2\lambda}
\int_0^1{du\over 1-u\mu}(1+{1\over 2\lambda}\log{1\over u})^m.
\ee
We break it into two regions of integration $0< u<{1\over 2}$ and
${1\over 2}\leq u\leq 1$. In the first region, the integral is of size
\bea
\int_0^{1\over 2}{du\over 1-u\mu}(1+{1\over 2\lambda}\log{1\over u})^m
&\leq&
2\int_0^{1\over 2}(1+{1\over 2\lambda}\log{1\over u})^m
\nonumber\\
&\leq&
{2\over\lambda^m}\int_0^{e^{-\lambda}}(\log{1\over u})^m du
+2^{m+1}\int_{e^{-\lambda}}^{1\over 2}du\leq C_m.
\eea
In the second region, we have
\be
\int_{1\over 2}^1{du\over 1-u\mu}(1+{1\over 2\lambda}\log{1\over u})^m
\leq
C_m
\int_{1\over 2}^1{du\over 1-u\mu}
\ee
This last integral can be evaluated explicitly, and we find that
it is bounded by $(1+\log (1+{1\over \lambda y})$.
This is the desired estimate.

\medskip
Next, consider the contribution of the region $t<-y$.
In this region, we have instead
\be
{1\over 2}
e^{\lambda t}
\leq {e^{\lambda y}+e^{-\lambda y}\over e^{\lambda(y- t)}+e^{-\lambda(y- t)}}
\leq 2{e^{\lambda y}\over e^{\lambda(y-t)}}\leq 2 e^{\lambda t}
\ee
The contribution from $t<-y$ to the left hand side of the desired inequality
can then be bounded by
\be
||\phi||_{\Lambda_{(m)}^0}(1+y)^m
\int_{-\infty}^{-y}{e^{2\lambda t}\over 1-e^{2\lambda t}}(1+|t|)^m\,dt
\leq C_m\,
{1\over\lambda}(1+\log (1+{1\over\lambda y}),
\ee
as was to be shown. Q.E.D.

\bigskip
The bound provided by Lemma \ref{lemmatinfinity} involves a $\log\,(\lambda y)^{-1}$ term,
and is not adequate
for $y$ close to $0$. This is because the integral is only a principal
value integral when $|t-y|$ is small, and the estimates we have just
derived for the contribution of the region $t>y$ do not take into
account the cancellations inherent to principal value integrals.
This issue is addressed in the next lemma:

\begin{lemma}
\label{lemmay0}
Assume that $0<\lambda y<1$. Then
\bea
&&
\bigg|\int_{y<|t|<1}{1\over e^{\lambda t}-e^{-\lambda t}}
{e^{\lambda y}+e^{-\lambda y}\over
e^{\lambda(y-t)}+e^{-\lambda(y-t)}}
\phi(y-t)dt
\bigg|
\leq C\,||\phi||_{C^0_{[0,2]}}(y+{1\over\lambda})
\\
&&
\bigg|\int_{|t|>1}{1\over e^{\lambda t}-e^{-\lambda t}}
{e^{\lambda y}+e^{-\lambda y}\over
e^{\lambda(y-t)}+e^{-\lambda(y-t)}}
\phi(y-t)dt
\bigg|
\leq C_m\,||\phi||_{\Lambda_{(m)}^0}\,{1\over\lambda}.
\eea
\end{lemma}

\noindent
{\it Proof.} Since we can assume that $\lambda$ is large,
the condition that $\lambda y<1$ implies that $y<1$, say.
We can exploit the cancellation by writing the integral
over the region $y<|t|<1$ in the form,
\be
\int_{y<|t|<1}
=
\int_y^1{1\over e^{\lambda t}-e^{-\lambda t}}
\{{e^{-\lambda y}+e^{\lambda y}\over e^{-\lambda(y-t)}+e^{\lambda(y-t)}}
\phi(y-t)
-
{e^{-\lambda y}+e^{\lambda y}\over e^{-\lambda(y+t)}+e^{\lambda(y+t)}}
\phi(y+t)\}dt
\ee
Next, the expression within brackets is written as,
\bea
&&
{e^{-\lambda y}+e^{\lambda y}\over e^{-\lambda(y-t)}+e^{\lambda(y-t)}}
\phi(y-t)
-
{e^{-\lambda y}+e^{\lambda y}\over e^{-\lambda(y+t)}+e^{\lambda(y+t)}}
\nonumber\\
&&
\qquad
={e^{-\lambda y}+e^{\lambda y}\over e^{-\lambda(y-t)}+e^{\lambda(y-t)}}
(\ \phi(y-t)-\phi(y+t)\ )
\nonumber\\
&&
\qquad
\qquad
+
\phi(y+t)
\{
{e^{-\lambda y}+e^{\lambda y}\over e^{-\lambda(y-t)}+e^{\lambda(y-t)}}
-
{e^{-\lambda y}+e^{\lambda y}\over e^{-\lambda(y+t)}+e^{\lambda(y+t)}}
\}
\eea
The contribution of the first term on the right hand side can be estimated
as follows,
\bea
\bigg|\int_y^1
{1\over e^{\lambda t}-e^{-\lambda t}}
{e^{-\lambda y}+e^{\lambda y}\over e^{-\lambda(y-t)}+e^{\lambda(y-t)}}
(\ \phi(y-t)-\phi(y+t)\ )
\bigg|
\leq
||\phi||_{C^1_{[-2,2]}}
\int_y^1 {t\over e^{\lambda t}-e^{-\lambda t}}e^{\lambda(2y-t)}dt
\nonumber
\eea
Since $\lambda y<1$, we can estimate this last term crudely by
\be
\int_y^1 {t\over e^{\lambda t}-e^{-\lambda t}}e^{\lambda(2y-t)}dt
\leq
{e^2\over\lambda}\int_{y}^1{t\lambda\over e^{2\lambda t}-1}dt
\leq
{C\over \lambda},
\ee
since the function $u(e^{2u}-1)^{-1}$ is a smooth and bounded function for $u\geq 0$.
Next, to estimate the other contribution,
we also exhibit the cancellation more clearly,
\bea
\label{group}
{1\over e^{-\lambda(y-t)}+e^{\lambda(y-t)}}
-
{1\over e^{-\lambda(y+t)}+e^{\lambda(y+t)}}
&=&
{1\over 1+e^{-2\lambda(t+y)}}(\ {1\over e^{\lambda(t-y)}}-{1\over e^{\lambda(t+y)}})
\nonumber\\
&&
\qquad
+
{1\over e^{\lambda(t+y)}}({1\over 1+e^{-2\lambda(t-y)}}
-
{1\over 1+e^{-2\lambda(y+t)}})
\nonumber
\eea
The first resulting group of terms can be estimated by
\be
\bigg|
{1\over 1+e^{-2\lambda(t+y)}}(\ {1\over e^{\lambda(t-y)}}-{1\over e^{\lambda(t+y)}})
\bigg|
\leq
e^{-\lambda t}(e^{\lambda y}-e^{-\lambda y})
\leq C\, \lambda y e^{-\lambda t},
\ee
and the corresponding integral in turn by,
\be
\int_y^1
dt\,{1\over e^{\lambda t}-e^{-\lambda t}}
|\phi(y+t)|
\cdot
\bigg|
{1\over 1+e^{-2\lambda(t+y)}}(\ {1\over e^{\lambda(t-y)}}-{1\over e^{\lambda(t+y)}})
\bigg|
\leq
||\phi||_{C^0_{[0,2]}}\,\lambda y\,
\int_y^1{dt\over e^{2\lambda t}-1}.
\ee
To determine the size of this expression, we break it up as follows,
\bea
\int_y^1{dt\over e^{2\lambda t}-1}
&=&
\int_y^{1\over \lambda}{dt\over e^{2\lambda t}-1}
+
\int_{1\over\lambda}^1{dt\over e^{2\lambda t}-1}
\leq
C\,(\ \int_y^{1\over\lambda}{dt\over\lambda t}
+
\int_{1\over\lambda}^1{dt\over e^{2\lambda t}}\ )
\nonumber\\
&\leq&
C\, {1\over\lambda}\,(\log{1\over\lambda y}+1),
\eea
and hence, since $\lambda y<1$,
\be
\int_y^1
dt\,{1\over e^{\lambda t}-e^{-\lambda t}}
|\phi(y+t)|
\cdot
\bigg|
{1\over 1+e^{-2\lambda(t+y)}}(\ {1\over e^{\lambda(t-y)}}-{1\over e^{\lambda(t+y)}})
\bigg|
\leq
||\phi||_{C^0_{[0,2]}}\,({1\over\lambda}+y).
\ee
The remaining group of terms in (\ref{group}) can be estimated in a similar way,
\be
{1\over e^{\lambda(t+y)}}
\bigg|{1\over 1+e^{-2\lambda(t-y)}}
-
{1\over 1+e^{-2\lambda(y+t)}}
\bigg|
\leq (e^{\lambda y}-e^{-\lambda y})e^{-\lambda(3t+y)}
\leq C\,\lambda y e^{-3\lambda t},
\ee
and hence
\bea
&&
\int_y^1
dt\,{1\over e^{\lambda t}-e^{-\lambda t}}
{1\over e^{\lambda(t+y)}}
\bigg|{1\over 1+e^{-2\lambda(t-y)}}
-
{1\over 1+e^{-2\lambda(y+t)}}
\bigg|\cdot
|\phi(y+t)|
\nonumber\\
&&
\quad
\leq
||\phi||_{C^0_{[0,2]}}
\lambda y\int_y^1{e^{-3\lambda t}\over e^{\lambda t}-e^{-\lambda t}}dt,
\eea
which is even smaller than the previous integral. Finally,
to estimate the integral from the region $|t|>1$, we
have the simple estimate, since $\lambda y<1$, say for $t>0$,
\bea
{1\over e^{\lambda t}-e^{-\lambda t}}
{e^{\lambda y}+e^{-\lambda y}
\over
e^{-\lambda(t-y)}+e^{\lambda(t-y)}}|\phi(y-t)|
&\leq&
C\, {1\over e^{\lambda t}}{1\over e^{\lambda(t-y)}}
||\phi||_{\Lambda_{(m)}^0}(1+|t|)^m
\nonumber\\
&\leq& C\,
||\phi||_{\Lambda_{(m)}^0}(1+|t|)^m e^{-2\lambda t}
\eea
which implies readily the desired inequality upon integration in $t$.
The proof of the lemma is complete. Q.E.D.

\bigskip
\noindent
{\it Proof of Theorem 2}. It suffices to combine all estimates
from Lemmas 4,5, and 6: when
$\lambda y\geq 1$, we apply Lemmas 4 and 5, and when $\lambda y<1$, we apply
Lemma 4 and 6. Q.E.D.


\newpage
\begin{thebibliography}{99}

\bibitem{B} Birkhoff, G.D.,
``The generalized Riemann problem for linear differential
equations and allied problems for linear difference and
q-difference equations",
Proc. Amer. Acad. of Arts and Sciences {\bf 49} (1913) 521-568.

\bibitem{Bo} Borodin, A.,
``Isomonodromy transformations of linear systems of
difference equations", Ann. of Math. {\bf 160} (2004) 1141-1182.

\bibitem{HI} Harnad, J. and A.R. Its,
``Integrable Fredholm operators and dual isomonodromic deformations",
Comm. Math. Phys. {\bf 226} (2002) 497-530.

\bibitem{K} Krichever, I.M.,
``Analytic theory of difference equations with rational and
elliptic coefficients and the Riemann-Hilbert problem",
Russian Math. Surveys {\bf 59} (2004) 1117-1154.

\bibitem{LU}
Ladyzhenskaya, O. and N. Uraltseva,
``{\it Linear and quasilinear elliptic equations}",
Academic Press, 1968.

\bibitem{Stein} Stein, E.M.,
``{\it Harmonic Analysis: real-variable methods, orthogonality,
and oscillatory integrals}",
Princeton University Press, 1993.

\bibitem{TW} Tracy, C.A. and H. Widom,
``Fredholm determinants, differential equations,
and matrix models", Comm. Math. Phys. {\bf 163} (1994) 33-72.


\end{thebibliography}
\end{document}